\input opmac
\input ams-math

\typosize[11/13.6]

\hyperlinks{\Blue}{\Green}

\def\hlink[#1]#2{{\def\tmpb{#2}%
   \replacestrings{//}{{\urlskip\urlslashslash\urlbskip}}%
   \replacestrings{/}{{\urlskip/\urlbskip}}%
   \replacestrings{.}{{\urlskip.\urlbskip}}%
   \replacestrings{?}{{\urlskip?\urlbskip}}%
   \replacestrings{=}{{\urlskip=\urlbskip}}%
   \replacestrings{~}{{\char`\~}}%
   \replacestrings{_}{{\char`\_}}%
   \replacestrings{^}{{\char`\^}}%
   \replacestrings{\\}{\bslash}%
   \replacestrings{\{}{{\char`\{}}%
   \replacestrings{\}}{{\char`\}}}%
   \replacestrings{&}{{\urlbskip\char`\& \urlskip}}%
   \def\|{}\ulink[#1]{\tmpb}%
}}

\eoldef\sec#1{\ifnonum\else \global\advance\secnum by1 \fi
  \sechook {\globaldefs=1 \seccnum=0}\relax
  \edef\thesecnum{\othe\chapnum.\the\secnum} 
  \def\dotocnumafter{\wtotoc1\rm{#1}}%
  \printsec{#1}\resetnonumnotoc
}

\def\so{{\rm so}}
\def\AW{{\rm AW}}

\def\tr{\mathop{\rm Tr}}

\def\im{{\rm i}}

\def\P{{\bbchar P}}
\def\N{{\bbchar N}}

\def\C{{\bbchar C}}
\def\Z{{\bbchar Z}}
\def\X{{\script X}}
\def\Y{{\script Y}}

\newcount\numA \newcount\numB \newcount\numC \newcount \numD

\def\numberedpar#1#2#3{\par\removelastskip\smallbreak \global\advance\csname num#1\endcsname by1
   \noindent\wlabel{\the\csname num#1\endcsname}%
   {\bf#2 \the\csname num#1\endcsname}\ifx^#3^\else{} (#3)\fi{\bf .}\space} 
\def\define{\numberedpar A{Definition}{}}
\def\theorem{\numberedpar B{Theorem}{}}
\def\lemma{\numberedpar C{Lemma}{}}
\def\nazlemma(#1) {\numberedpar C{Lemma}{#1}}
\def\remark{\numberedpar D{Remark}{}}
\def\proof{\par\noindent{\bf Proof.}\space }
\def\qed{\leavevmode\hbox{\quad}\hfill$\square$\smallbreak}

\let\phi=\varphi
\let\epsilon=\varepsilon
\let\theta=\vartheta
\let\rho=\varrho

\def\Fi#1#2(#3;#4|#5;#6){{}_{#1}\phi_{#2}\left(\matrix{#3 \cr #4}\biggm|{#5};{#6}\right)}

\tit Representations of Askey--Wilson algebra

\centerline{Daniel Gromada$^{\dag}$, Severin Pošta$^{\ddag}$}
\bigskip

{\parindent=0pt
$^{\dag}$~Department of Physics, Faculty of Nuclear Sciences and Physical Engineering,

$\phantom{^{\dag}}$~Czech Technical University in Prague, B\v rehov\'a 7, CZ-115 19 Prague, Czech Republic

$\phantom{^{\dag}}$~E-mail: d.gromada@seznam.cz

$^{\ddag}$~Department of Mathematics, Faculty of Nuclear Sciences and Physical Engineering,

$\phantom{^{\ddag}}$~Czech Technical University in Prague, Trojanova 13, CZ-120 00 Prague, Czech Republic

$\phantom{^{\ddag}}$~E-mail: severin.posta@fjfi.cvut.cz
}

\bigskip
\noindent\hfil\vbox{\hsize=0.8\hsize \parindent=0pt
{\bf Abstract.}\space We deal with the classification problem of finite-dimensional representations of so called Askey--Wilson algebra in the case when $q$ is not a root of unity. We classify all representations satisfying certain property, which ensures diagonalizability of one of the generating elements.

}

\sec Introduction

The relationship between special classes of orthogonal polynomials from the Askey scheme is already known for some time. Such a relationship is very valuable since the basic properties of the orthogonal polynomials can be derived from this relationship very easily. In 1991 \cite[Zhedanov] Alexei Zhedanov constructed a $q$-commutator algebra that correspond to the most general class of discrete polynomials from the Askey--Wilson scheme---the $q$-Racah polynomials defined by Askey and Wilson in \cite[AskeyWilson] (see \cite[Koekoek] for detailed overview of the Askey--Wilson scheme). Zhedanov named the algebra after authors of the scheme Richard Askey and James Wilson.

The important concept that emerges in both Askey--Wilson algebra and finite orthogonal polynomial sequences and connects those structures is so called Leonard pair, which is a pair of operators such that both are tridiagonal in the eigenbasis of the other. This correspondence was first discovered by Leonard \cite[Leonard]. The pair of such operators was named after him by Terwilliger \cite[TerwilligerLeonard]. For a nice introduction to the theory of Leonard pairs, see \cite[TerwilligerIntro]. The Leonard pair corresponding to $q$-Racah polynomials is studied in \cite[TerwilligerRacah].

In the Askey--Wilson algebra the Leonard pair is made of elements of the algebra in certain representation. If we want to find the Leonard pair and discover the orthogonal polynomials and their properties without using the properties of the orthogonal polynomials in the first place, we can make use of classification of representations and certain automorphisms, as we have shown for example in \cite[GromadaPosta]. Although Zhedanov was able to construct suitable representation of the Askey--Wilson algebra and show the connection with $q$-Racah polynomials, the complete classification was not available until recent work of Hau-Wen Huang \cite[Huang]. (Only in the case when $q$ is not a root of unity. The case when $q$ is a root of unity is much more complicated and the classification problem is still open, see \cite[Uqso3,centeruqso3,centeraw].)

Independently on Huang we worked on the classification problem for Askey--Wilson algebra as well. Our approach follows the Zhedanov's paper \cite[Zhedanov] using the technique of shift operators to construct the representations. However, it completes the Zhedanov's paper with rigorous mathematical theorems and proofs. The technique is inspired by classification of representations of algebra $U'_q(\so_3)$ \cite[Uqso399,Uqso3], which is a special case of Askey--Wilson algebra.

We show that every representation satisfying certain property (certain numbers should not be eigenvalues of one of the generating element) can be constructed using shift operators and that those representations are determined by trace of the generating element (and parameters of the algebra of course). Consequently, we were able to prove that the dual representation constructed by Zhedanov, which provides the relationship with $q$-Racah polynomials, is indeed equivalent to the original one without need of prior knowledge of the properties of the $q$-Racah polynomials.

We also consider separately special case when certain parameters of the Askey--Wilson algebra are zero, which allows to construct new representations by the shift operators not considered by Zhedanov.

Although complete classification of Askey--Wilson algebra representations was already presented, we think that our approach is still interesting since it is very straightforward using the classical technique of shift operators and extends the Zhedanov's paper.

In the whole work, we assume that $q$ is not a root of unity.

\sec Definition and basic properties

The definition of Askey--Wilson algebra was presented by Alexey Zhedanov in \cite[Zhedanov]. We will use a $\Z_3$-symmetric presentation, which was mentioned for example in \cite[Wiegmann]. We chose the presentation in a way that the algebra $U'_q(\so_3)$ is a special case when all parameters equal to zero.

\label[D.AW]
\define Askey--Wilson algebra $\AW_q(A_1,A_2,A_3)$ is a complex associative algebra generated by three elements $I_1$, $I_2$, $I_3$ and relations
$$\eqalignno{
\label[eq.AW1]	q^{1/2}I_1I_2-q^{-1/2}I_2I_1&=I_3+A_3,&\eqmark\cr
\label[eq.AW2]	q^{1/2}I_2I_3-q^{-1/2}I_3I_2&=I_1+A_1,&\eqmark\cr
\label[eq.AW3]	q^{1/2}I_3I_1-q^{-1/2}I_1I_3&=I_2+A_2,&\eqmark\cr
}$$
where $A_1$, $A_2$ and $A_3$ are complex parameters. We will often denote the algebra only shortly as $\AW$.

Later, Paul Terwilliger in \cite[TerwilligerAW] presented so called {\em universal} Askey--Wilson algebra, that was also defined by three generating elements satisfying the same commutation relations, but $A_1$, $A_2$ and $A_3$ are not considered as complex parameters, but as elements of the center of the algebra.


\theorem The set of monomials $\{I_1^kI_2^mI_3^n\mid k,m,n\in{\bbchar N}_0\}$ forms a basis of $\AW$.
\proof As in the case of $U'_q(\so_3)$ \cite[Uqso3] or universal Askey--Wilson algebra \cite[TerwilligerAW] we make use of the Diamond lemma \cite[Bergman1978]. We transform the generating relation into a form compatible with ordering $I_1\le I_2\le I_3$:
$$\eqalign{I_2I_1&=qI_1I_2-q^{1/2}(I_3+A_3),\cr
I_3I_2&=qI_2I_3-q^{1/2}(I_1+A_1),\cr
I_3I_1&=q^{-1}I_1I_3+q^{-1/2}(I_2+A_2).}$$
We see that there are no inclusion ambiguities on the left hand side and there is only one overlap ambiguity---a monomial $I_3I_2I_1$ can be reduced using the first or the second relation. We show that this ambiguity is resolvable.

Reducing by the first relation we get
$$\eqalignno{I_3I_2I_1&=qI_3I_1I_2-q^{1/2}I_3(I_3+A_3)=I_1I_3I_2+q^{1/2}(I_2(I_2+A_2)-I_3(I_3+A_3))\cr
&=qI_1I_2I_3+q^{1/2}(-I_1(I_1+A_1)+I_2(I_2+A_2)-I_3(I_3+A_3)),\cr
\noalign{\hbox{while the second relation leads to}}
&=qI_2I_3I_1-q^{1/2}I_3(I_1+A_1)=I_2I_1I_3+q^{1/2}(-I_1(I_1+A_1)+I_2(I_2+A_2))\cr
&=qI_1I_2I_3+q^{1/2}(-I_1(I_1+A_1)+I_2(I_2+A_2)-I_3(I_3+A_3)).}$$
Both results are the same, so the ambiguity is resolvable.

Finally, we can see that $\{I_1^kI_2^mI_3^n\}$ is indeed the set of all reduced monomials.\qed

\label[L.AWiso]
\lemma There are the following isomorphisms of Askey--Wilson algebras
$$\rho\colon\AW_q(A_1,A_2,A_3)\to\AW_q(A_2,A_3,A_1),\quad I_1\mapsto I_2,\quad I_2\mapsto I_3,\quad I_3\mapsto I_1,$$
$$\sigma\colon\AW_q(A_1,A_2,A_3)\to\AW_q(A_2,A_1,A_3),$$
$$\quad I_1\mapsto I_2,\quad I_2\mapsto I_1,\quad I_3\mapsto I_3+(I_2I_1-I_1I_2)(q^{1/2}+q^{-1/2}),$$
$$\tau_{\epsilon,\epsilon'}\colon\AW_q(A_1,A_2,A_3)\to\AW_q(\epsilon A_1,\epsilon' A_2,\epsilon\epsilon' A_3),\quad I_1\mapsto\epsilon I_1,\quad I_2\mapsto\epsilon' I_2,\quad I_3\mapsto\epsilon\epsilon' I_3,$$
where $\epsilon, \epsilon'\in\{1,-1\}$.

\smallskip
Those isomorphisms can be also interpreted as automorphisms of the universal Askey--Wilson algebra. In fact, the first two form a faithful action of the group $\rm PSL_2({\bbchar Z})$ \cite[TerwilligerAW].

We are going to use the following definition of {\em $q$-numbers}
$$[\alpha]_q:={q^\alpha-q^{-\alpha}\over q-q^{-1}}\quad\hbox{for $\alpha\in\C$}.\eqmark$$

\label[L.qnum]
\nazlemma(Lemma 10 \cite[Uqso3]) Let $q$ is not a root of unity.
\begitems \style n
	* If $q^{2\alpha}=-q^l$ for some $l\in{\bbchar Z}$, then $[\alpha-j]_q=[\alpha-k]_q$ for $j,k\in{\bbchar Z}$, $j\neq k$, if and only if $j+k=l$.
	* Conversely, having $[\alpha-j]_q=[\alpha-k]_q$ for some $\alpha\in\C$ a~$j,k\in{\bbchar Z}$, $j\neq k$, it follows that $q^{2\alpha}=-q^{j+k}$.
	* For every $\lambda\in\C$ there exists $\alpha\in\C$ such that $\lambda=[\alpha]_q$ and numbers $[\alpha]_q, [\alpha+1]_q,[\alpha+2]_q,\dots$ are mutually different.
\enditems

\sec Classification of finite-dimensional representations

As we indicated in the introduction, we follow the Zhedanov's construction of the representations. So, the definitions of Casimir element, shift operators, characteristic polynomial or dual representations are all inspired by the original paper \cite[Zhedanov].

\label[L.casimir]
\lemma {\em Casimir element}
\label[eq.casimir]
$$C:=q^2I_1^2+I_2^2+q^2I_3^2-(q^{5/2}-q^{1/2})I_1I_2I_3+q(q+1)A_1I_1+(q+1)A_2I_2+q(q+1)A_3I_3\eqmark$$
is a central element of $\AW$.

\define For suitable complex numbers $\lambda$, let $O_\lambda$ and $R_\lambda$ be the following linear combinations:
$$\eqalignno{
\label[eq.Odef]	O_\lambda&:={-\im q^{1/2}A_1+q^{-\lambda}A_2\over [\lambda]_q([\lambda]_q-[\lambda+1]_q)}I_3+\im I_2+q^{-\lambda+1/2}I_1,&\eqmark\cr
\label[eq.Rdef]	R_\lambda&:={-\im q^{1/2}A_1-q^{\lambda}A_2\over [\lambda]_q([\lambda]_q-[\lambda-1]_q)}I_3+\im I_2-q^{\lambda+1/2}I_1.&\eqmark
}$$
In case of $A_1=A_2=0$ we consider this definition without the first term and we assume that $\lambda$ is an arbitrary complex number. Otherwise, we require the expression to be well defined, so $[\lambda]_q\neq[\lambda+1]_q$ and $[\lambda]_q\neq[\lambda-1]_q$, respectively.
\smallskip

In the following text, considering a representation $R$ of the algebra $\AW$ on a vector space $V$, we will denote the representing linear operators just $I_1$, $I_2$, and $I_3$ instead of $R(I_1)$, $R(I_2)$, and $R(I_3)$.

\lemma Let $R$ be a representation of $\AW$ on $V$. Let $x\in\ker(I_3+\im [\lambda]_q)$. Then
$$\eqalignno{
\label[eq.Oshift]	I_3(O_\lambda x)&=-\im [\lambda+1]_qO_\lambda x,&\eqmark\cr
\label[eq.Rshift]	I_3(R_\lambda x)&=-\im [\lambda-1]_qR_\lambda x,&\eqmark\cr
\label[eq.OR]	O_{\lambda-1}R_\lambda x&=(\tilde C_{\lambda-1}-C)x,&\eqmark\cr
\label[eq.RO]	R_{\lambda+1}O_\lambda x&=(\tilde C_\lambda-C)x,&\eqmark\cr
}$$
where
$$\tilde C_\lambda=-q[\lambda]_q[\lambda+1]_q-\im q([\lambda]_q+[\lambda+1]_q)A_3-q{A_1^2+\im (q^{-\lambda-1/2}-q^{\lambda+1/2})A_1A_2+A_2^2\over([\lambda]_q-[\lambda+1]_q)^2}.$$
\proof By inspection using the $q$-commutation relations \ref[eq.AW1]--\ref[eq.AW3].\qed

\label[L.I1I2soust]
\lemma Take $y,z\in V$ and $\lambda\in\C$ such that both $O_\lambda$ and $R_\lambda$ are well-defined. Then the system of linear equations
$$O_\lambda x=y,\quad R_\lambda x=z$$
for vectors $I_1x$ and $I_2x$ has a unique solution if and only if $q^\lambda\neq\im\epsilon$ for $\epsilon\in\{-1,1\}$.

\proof From \ref[eq.Odef] and \ref[eq.Rdef] we can easily express $(q^{-\lambda}+q^{\lambda})I_1x$. The factor in the bracket is nonzero if and only if $q^\lambda\neq\im\epsilon$. The same holds for $I_2x$.\qed

Now, consider a general representation of Askey--Wilson algebra. We construct a basis in which $I_3$ acts diagonally and we try to express the explicit form of such representation. The construction of the basis is an analogy of Theorem 3 from \cite[Uqso3], where the same procedure is performed for $U'_q(\so_3)$. To perform such a construction, we need to assume that $2\epsilon/(q-q^{-1})$ and $\epsilon/(q^{1/2}-q^{-1/2})$, $\epsilon=\pm 1$ are not an eigenvalues of $I_3$. Such representations will be called {\em classical} (inspired by notation in \cite[Uqso3]).

These conditions allow us to use the shift operators to construct the eigenbasis of $I_3$. The first condition ensures the existence of solution of the system in lemma \ref[L.I1I2soust] and the second one ensures the shift operators to be well defined.
However, in the case when $A_1=A_2=0$, which will be considered separately, we will be able to construct the eigenbasis of $I_3$ even for the non-classical representations.



\label[secc.obrep]
\secc Classical representations

\label[V.AWbaze]
\theorem Let $q$ is not a root of unity. Let $R$ be an irreducible classical representation (i.e. satisfying  $\ker(I_3-2\epsilon/(q-q^{-1}))=\{0\}$ and $\ker(I_3-\epsilon/(q^{1/2}-q^{-1/2}))\neq\{0\}$ for $\epsilon=\pm 1$) of $\AW_q(A_1,A_2,A_3)$ on $V$, $\dim V=N+1$. Then there exists a complex number $\mu$ and a non-zero vector $v_0\in\ker(I_3+\im[\mu]_q)$ such that $O_\mu v_0=0$. If we define
$$v_{j+1}:=R_{\mu-j}v_j\quad\hbox{for $j=0,1,\dots r-1$},$$
then the tuple $(v_0,\dots,v_N)$ forms a basis of $V$.
\proof Let $w_0$ be an eigenvector corresponding to eigenvalue $-\im[\tilde\mu]_q$, where $\tilde \mu$ is chosen in such a way that the numbers $-\im[\tilde\mu]_q, -\im[\tilde\mu+1]_q,\dots$ are mutually different (see Lemma \ref[L.qnum]). Now, define
$$w_{j+1}:=O_{\tilde\mu+j}w_j\quad\hbox{for $j\ge 0$}$$
and denote $l\in\N$ such that $w_0,\dots,w_{l-1}$ are linearly independent, but $w_l$ is already linearly dependent on them. Such $l$ has to exist since $V$ is finite-dimensional. From the equation \ref[eq.Oshift] it follows that the vectors $w_j$ satisfy the eigenequation $I_3w_j=-\im[\tilde\mu+j]_qw_j$, where all the eigenvalues are mutually different. Since eigenvectors corresponding to mutually different eigenvalues are linearly independent, we have $w_l=0$.

Denote $v_0:=w_{l-1}$ a~$\mu:=\tilde\mu+l-1$ and define the following vectors:
\label[eq.vjdef]
$$v_{j+1}:=R_{\mu-j}v_j\quad\hbox{for $j\ge 0$.}\eqmark$$
From the equation \ref[eq.Rshift] it follows that they also have to satisfy the eigenequations
\label[eq.I3v]
$$I_3v_j=-\im[\mu-j]_qv_j.\eqmark$$
The possibility of $[\mu-j]_q=[\mu-j-1]_q$ for some $j$, so $R_{\mu-j}$ would not be well defined, contradicts the assumptions since it would follow that $q^\mu=\im\epsilon q^{j-1/2}$, so $-\im[\mu-j]_q=\epsilon/(q^{1/2}-q^{-1/2})$ would be an eigenvalue of $I_3$. Now we can denote $k\in\N$ such that the tuple $v_0,\dots,v_{k-1}$ is linearly independent, while $v_k$ is already linearly dependent. From \ref[eq.OR] it follows that
\label[eq.Ov]
$$O_{\mu-j}v_j=O_{\mu-j}R_{\mu-j+1}v_{j-1}=(\tilde C_{\mu-j}-C)v_{j-1},\quad O_\mu v_0=y_l=0\eqmark$$
for $j=1,2,\dots,k-1$.

From the irreducibility of the representation it follows that $C$ is a multiple of identity. This complex number will be also denoted by $C$. Its value is determined by equality $0=R_{\mu+1}O_\mu v_0=(\tilde C_\mu-C)v_0$, so $C=\tilde C_\mu$.

The equations \ref[eq.vjdef] and \ref[eq.Ov] define a system of equations
\label[eq.I1I2soust]
$$R_{\mu-j}v_j=v_{j+1},\quad O_{\mu-j}v_j=(\tilde C_{\mu-j}-C)v_{j-1},\eqmark$$
which, according to Lemma \ref[L.I1I2soust], has a solution for $I_1v_j$ a~$I_2v_j$ if and only if $q^{\mu-j}\neq\im\epsilon$ for $\epsilon\in\{-1,1\}$. This condition is satisfied thanks to our assumptions, otherwise $[\mu-j]_q=2\epsilon/(q-q^{-1})$ would be an eigenvalue of $I_3$ corresponding to the eigenvector $v_j$. Therefore, we were able to express the vectors $I_1v_j$ and $I_2v_j$ as a linear combination of $v_0,\dots,v_{k-1}$. Moreover, $I_3$ acts diagonally, so the span of $\{v_0,\dots,v_{k-1}\}$ is an invariant subspace and thanks to the irreducibility it has to be equal to the whole space $V$.
\qed

\label[L.klasvlč]
\lemma For arbitrary $(N+1)$-dimensional classical irreducible representation the operator $I_3$ is diagonalizable and has mutually different eigenvalues. Its eigenvalues are $-\im[\mu]_q$, $-\im[\mu-1]_q$, \dots, $-\im[\mu-N]_q$. Denoting $v_0$, \dots, $v_N$ the corresponding eigenvectors, we have $O_\mu v_0=0$ and $v_r:=R_{\mu-N}v_N=0$. In addition, it holds that $q^{2\mu}\neq q^l$ for every $l\in\{-1,0,1,\dots,2N+1\}$.
\proof The diagonalizability and spectrum of $I_3$ follows from the previous theorem. From assumptions, we have $-\im[\mu-j]_q\neq 2\epsilon/(q-q^{-1})$, so $q^\mu\neq\im\epsilon q^j$ for $j=0,\dots,N$. We also assume that $-\im[\mu-j]_q\neq\epsilon/(q^{1/2}-q^{-1/2})$, which means $q^\mu\neq\im\epsilon q^{j\pm 1/2}$ for $j=0,\dots,N$. Together it follows that $q^\mu\neq\im\epsilon q^{l/2}$, $l=-1,\dots,2N+1$.

This condition implies that the numbers $[\mu+1]_q$, \dots, $[\mu-N-1]_q$ are mutually different according to Lemma \ref[L.qnum].

The vector $O_\mu v_0$ is either zero or an eigenvector corresponding to eigenvalue $-\im[\mu+1]_q$. The second possibility cannot take place since $-\im[\mu+1]_q$ is not in the spectrum of $I_3$. The same holds for the vector $R_{\mu-N}v_N$.\qed

In the following text, we will still consider a representation satisfying the assumptions of Theorem \ref[V.AWbaze].

We will often deal with sums and products of consecutive $q$-numbers. Hence, it will be useful to summarize following relations
$$([\lambda+1]_q\pm[\lambda]_q)(q^{1/2}\mp q^{-1/2})=q^{\lambda+1/2}\mp q^{-\lambda-1/2},$$
$$[\lambda]_q[\lambda+1]_q={([\lambda]_q+[\lambda+1]_q)^2-1\over (q^{1/2}+q^{-1/2})^2},$$
$$([\lambda]_q-[\lambda+1]_q)^2=([\lambda]_q+[\lambda+1]_q)^2{(q^{1/2}-q^{-1/2})^2\over(q^{1/2}+q^{-1/2})^2}+{4\over(q^{1/2}+q^{-1/2})^2}.$$

We denote $\tilde D_j:=\tilde C_\mu-\tilde C_{\mu-j}$ which will play an important role for characterization of the representations. Next, we denote
\label[eq.lambdadef]
$$\Lambda_j:=[\mu-j+1]_q+[\mu-j]_q={q^{\mu-j+1/2}-q^{-\mu+j-1/2}\over(q^{1/2}-q^{-1/2})},\eqmark$$
$$g_j:=[\mu-j+1]_q-[\mu-j]_q={q^{\mu-j+1/2}+q^{-\mu+j-1/2}\over(q^{1/2}+q^{-1/2})}.\eqmark$$
Using the relations above and the definition of $\tilde C_\mu$ we see that $\tilde D_jg_j^2$ is a polynomial of fourth order in variable $\Lambda_j$. The roots of this polynomial determine the form of the representation. Thus, we could call it a {\em characteristic polynomial}. (Characteristic polynomial is defined in \cite[Zhedanov] in a similar way.) From the definition of $\tilde D_j$ it is clear that $\tilde D_0=0$; thus, one of the roots of the characteristic polynomial is $\Lambda_0$.

Solving the system \ref[eq.I1I2soust] we get
\label[eq.I1v]
$$I_1v_j={-q^{-1/2}\over q^{-\mu+j}+q^{\mu-j}}(\tilde D_jv_{j-1}+v_{j+1})-{A_1-\im [\mu-j]_qA_2(q^{1/2}-q^{-1/2})\over g_jg_{j+1}}v_j,\eqmark$$
\label[eq.I2v]
$$I_2v_j={\im\over q^{-\mu+j}+q^{\mu-j}}(q^{\mu-j}\tilde D_jv_{j-1}-q^{-\mu+j}v_{j+1})-{-\im [\mu-j]_qA_1(q^{1/2}-q^{-1/2})+A_2\over g_jg_{j+1}}v_j.\eqmark$$

Denote $\Lambda_{j_0},\dots,\Lambda_{j_3}$ the roots of the characteristic polynomial and choose $j_0=0$. Since we know, what is the leading coefficient, we can factor the polynomial
\label[eq.Djjk]
$$\eqalign{\tilde D_j&={q(q^{1/2}-q^{-1/2})^2\over(q^{1/2}+q^{-1/2})^4}{\prod_{k=0}^3(\Lambda_j-\Lambda_{j_k})\over g_j^2}\cr
&=q{\prod_{k=0}^3(q^{\mu-j+1/2}-q^{-\mu+j-1/2}-q^{\mu-j_k+1/2}+q^{-\mu+j_k-1/2})\over(q-q^{-1})^2(q^{\mu-j+1/2}+q^{-\mu+j-1/2})^2}\cr
&={q^{2\mu-2j+2}\prod_{k=0}^3\bigl((1-q^{j-j_k})(1+q^{-2\mu+j+j_k-1})\bigr)\over(q-q^{-1})^2(1+q^{-2\mu+2j-1})^2}}\eqmark$$

\label[L.chpolkor]
\lemma For any $(N+1)$-dimensional classical irreducible representation the numbers $\Lambda_0$ and $\Lambda_r$ are mutually different roots of the characteristic polynomial. For any $k\in\{1,\dots,N\}$ the number $\Lambda_k$ is not a root of the characteristic polynomial.
\proof We can see that if $\tilde D_k=0$ for $k\in\{1,\dots,N\}$, it would mean that the span of $v_k,\dots,v_N$ is a non-trivial invariant subspace, which contradicts the irreducibility of the representations.

Substituting \ref[eq.I1v], \ref[eq.I2v] and \ref[eq.I3v] to the condition $(q^{1/2}I_1I_2-q^{-1/2}I_2I_1-I_3-A_3)v_N=0$ (from the relation \ref[eq.AW1]) we get the following relation
\label[eq.Dcond]
$${-\im q^{-1}\tilde D_{N+1}\over q^{\mu-N}+q^{-\mu+N}}v_N=0,\eqmark$$
which implies $\tilde D_{N+1}=0$, so $\Lambda_{N+1}$ is a root of the characteristic polynomial as well. Now, we have to show that it is not equal to $\Lambda_0$.

Assuming $\Lambda_0=\Lambda_{N+1}$ and using the fact that $q$ is not a root of unity we get $q^{2\mu}=-q^N$, which is excluded by Lemma \ref[L.klasvlč].
\qed

The form of the representation of course depends on the parameters of the algebra $A_1$, $A_2$, and $A_3$. Besides that, it also depends on the value of the Casimir element $C$, which is determined by the number $q^\mu$. The possible values of $q^\mu$ are restricted by the assumption of finite dimension $\tilde D_{N+1}=0$. This is a polynomial equation of degree eight in variable $q^\mu$. Excluding the possibility $\Lambda_{N+1}=\Lambda_0$, we can reduce the degree to six by eliminating the factor $\Lambda_{N+1}-\Lambda_0$. However, the equation remains very hard to solve and we will not try to express the possible values of $q^\mu$ explicitly.

Representations with different $q^\mu$ may be equivalent. We will show that for fixed parameters of the algebra, the class of $(N+1)$-dimensional equivalent representations is determined by number $[\mu-N/2]_q$, which is, up to a constant, the trace of $I_3$ (traces were used to distinguish representations already in the case of $U'_q(\so_3)$ \cite[Uqso399] and they are also used in the Huang's classification \cite[Huang]).

\label[L.Rekviv]
\lemma Let $R_1$ and $R_2$ be classical irreducible $r$-dimensional representations of $\AW$ constructed in Theorem \ref[V.AWbaze], $\{v_j^{(1)}\}$ and $\{v_j^{(2)}\}$ the corresponding bases, and $\mu_1$ and $\mu_2$ the corresponding complex numbers characterizing the representations. Then the representations $R_1$ and $R_2$ are equivalent if and only if $[\mu_1-N/2]_q=[\mu_2-N/2]_q$.
\proof
We compute the trace of $I_3$
\label[eq.tr]
$$\tr I_3=\sum_{j=0}^N-\im[\mu-j]_q=-\im{q^{(N+1)/2}-q^{-(N+1)/2}\over q^{1/2}-q^{-1/2}}[\mu-N/2]_q.\eqmark$$
Therefore, representations with different $[\mu-N/2]_q$ have to be inequivalent.

Conversely, let $[\mu_1-N/2]_q=[\mu_2-N/2]_q$. The equation \ref[eq.tr] is quadratic in $q^{\mu-N/2}$ so it can have two solutions. Generally, it holds that $[\lambda]_q=-[-\lambda]_q$. If $q^{\mu_1-N/2}$ is one of the solutions, the second one has to be $q^{\mu_2-N/2}=-q^{-\mu_1+N/2}$. This implies that
$$[\mu_1-j]_q=[\mu_2-N+j]_q,\quad q^{\mu_1-j}+q^{-\mu_1+j}=-(q^{\mu_2-N+j}+q^{-\mu_2+N-j}),\quad \Lambda^{(1)}_j=\Lambda^{(2)}_{N+1-j},$$
where $\Lambda^{(1)}_j$ a~$\Lambda^{(2)}_j$ correspond to the definition \ref[eq.lambdadef] for $\mu_1$ and $\mu_2$. Similarly, we denote $\tilde D_j^{(1)}$ and $\tilde D_j^{(2)}$ and we define a basis $\{w_0,\dots,w_{r-1}\}$ by equation
$$v_j^{(2)}=\prod_{k=1}^j\tilde D_k^{(1)}\cdot w_{N-j}.$$
Using \ref[eq.I1v], \ref[eq.I2v], and \ref[eq.I3v] we can check that the representation $R_1$ has the same matrix elements in the basis $\{v_j^{(1)}\}$ as $R_2$ in $\{w_j\}$.
\qed

To demonstrate the relationship of the representations with orthogonal polynomials it will be convenient to express the representation in terms of the numbers $j_1$, $j_2$, and $j_3$ instead of the parameters of the algebra $A_1$, $A_2$, and $A_3$. The parameters of the algebra together with the parameter of the representation $q^\mu$ determine the form of $\tilde D_j$, and hence the numbers $\Lambda_{j_1}$, $\Lambda_{j_2}$ a~$\Lambda_{j_3}$. 

Expressing the numbers $q^{j_1}$, $q^{j_2}$ a~$q^{j_3}$ in terms of the parameters of the algebra would lead to a bicubic equation, which would be very complicated and we will not perform it. To express the form of the representation in terms of those numbers, we will need the opposite relationship---to express the parameters of the algebra. This can be achieved using Vi\`et formulas for polynomial $\tilde D_j$ in the variable $\Lambda_j$. For given numbers $\mu$, $j_1$, $j_2$, $j_3$ we obtain the parameters $A_1$, $A_2$, $A_3$ such that $\Lambda_{j_1}$, $\Lambda_{j_2}$ a~$\Lambda_{j_3}$ are roots of $\tilde D_j$.

Using one of the relations, we are able to express
\label[eq.A3jk]
$$A_3={\im\over (q^{1/2}+q^{-1/2})^2}(\Lambda_0+\Lambda_{j_1}+\Lambda_{j_2}+\Lambda_{j_3}),\eqmark$$
and using the others we express a system of equations for $A_1$ and $A_2$
$$\eqalign{&(q^{1/2}-q^{-1/2})^3(\Lambda_0\Lambda_{j_1}\Lambda_{j_2}+\Lambda_{j_1}\Lambda_{j_2}\Lambda_{j_3}+\Lambda_{j_2}\Lambda_{j_3}\Lambda_0+\Lambda_{j_3}\Lambda_0\Lambda_{j_1})\cr-&4(q^{1/2}-q^{-1/2})(\Lambda_0+\Lambda_{j_1}+\Lambda_{j_2}+\Lambda_{j_3})=\im(q^{1/2}+q^{-1/2})^2(q-q^{-1})^2A_1A_2,}$$
$$\eqalign{-&4(q^{1/2}-q^{-1/2})^2\bigl(\Lambda_0(\Lambda_{j_1}+\Lambda_{j_2}+\Lambda_{j_3})+\Lambda_{j_1}(\Lambda_{j_2}+\Lambda_{j_3})+\Lambda_{j_2}\Lambda_{j_3}\bigr)\cr+&(q^{1/2}-q^{-1/2})^4\Lambda_0\Lambda_{j_1}\Lambda_{j_2}\Lambda_{j_3}+16=(q^{1/2}+q^{-1/2})^2(q-q^{-1})^2(A_1^2+A_2^2).}$$

This system leads to a biquadratic equation, so it has four solutions. From symmetry of the problem it is evident that if we denote $(A_1,A_2)$ one of the solutions, the other solutions are $(-A_1, -A_2)$, $(A_2,A_1)$ a~$(-A_2,-A_1)$. The explicit form of one of the solutions is following
$$\eqalignno{
\label[eq.A1jk]
A_1=&{\im\over(q^{1/2}+q^{-1/2})(q-q^{-1})}&\eqmark\cr
\bigl(&q^{(j_1+j_2+j_3-2\mu-1)/2}+q^{(-j_1+j_2+j_3-2\mu-1)/2}+q^{(j_1-j_2+j_3-2\mu-1)/2}+q^{(j_1+j_2-j_3-2\mu-1)/2}\cr
-&q^{(2\mu+1-j_1-j_2+j_3)/2}-q^{(2\mu+1-j_1-j_2-j_3)/2}-q^{(2\mu+1+j_1-j_2-j_3)/2}-q^{(2\mu+1-j_1+j_2-j_3)/2}\bigr),\cr
\label[eq.A2jk]
A_2=&{1\over(q^{1/2}+q^{-1/2})(q-q^{-1})}&\eqmark\cr
&\bigl(q^{-2\mu-1+(j_1+j_2+j_3)/2}-q^{(-j_1+j_2+j_3)/2}-q^{(j_1-j_2+j_3)/2}-q^{(j_1+j_2-j_3)/2}\cr
&-q^{-(-j_1+j_2+j_3)/2}-q^{-(j_1-j_2+j_3)/2}-q^{-(j_1+j_2-j_3)/2}+q^{2\mu+1-(j_1+j_2+j_3)/2}\bigr).}$$

We have found a representation for four isomorphic Askey--Wilson algebras with different parameters. If the necessary conditions for $(N+1)$-dimensional classical representation given by lemmata \ref[L.klasvlč] and \ref[L.chpolkor] are satisfied, then, by explicit computation, we can check that this representation really satisfies the commutation relations \ref[eq.AW1]--\ref[eq.AW3]. Since $I_3$ has mutually different eigenvalues, we can show the irreducibility of this representation easily. Taking an invariant subspace and an eigenvector of $I_3$ lying in this subspace. Using the shift operators $O_\lambda$ a~$R_\lambda$ we construct the remaining elements of the eigenbasis. Therefore, we have proven the following theorem.

\label[V.obecklas]
\theorem Let $q$ is not a root of unity. For arbitrary classical irreducible representation of Askey--Wilson algebra there exist complex numbers $\mu$, $j_1$, $j_2$, and $j_3$ such that this representation is equivalent to the representation given by equations \ref[eq.I3v], \ref[eq.I1v]--\ref[eq.Djjk]. Conversely, let the complex numbers $\mu$, $j_1$, $j_2$, and $j_3$ satisfy the following assumptions: $q^{2\mu}\neq- q^l$ for every $l\in\{-1,0,1,\dots,2N+1\}$, one of the numbers $\Lambda_{j_1}$, $\Lambda_{j_2}$, or $\Lambda_{j_3}$ equals to $\Lambda_{N+1}$, and none of those numbers equals to $\Lambda_k$ for all $k\in\{1,\dots,N\}$. Then the equations \ref[eq.I3v], \ref[eq.I1v]--\ref[eq.Djjk] define an irreducible representations of algebras $\AW_q(A_1,A_2,A_3)$, $\AW_q(-A_1,-A_2,A_3)$, $\AW_q(A_2,A_1,A_3)$, and $\AW_q(-A_2,-A_1,A_3)$, where the numbers $A_1$, $A_2$, and $A_3$ are determined by equations \ref[eq.A3jk]--\ref[eq.A2jk].

\label[P.assumptions]
\remark The assumptions can be somehow reformulated and a bit simplified. Firstly, the representation obviously does not depend on the order of the numbers $j_1$, $j_2$, $j_3$, so we can fix, for example $\Lambda_{j_3}=\Lambda_r$. Note, however, that the $q$-Racah polynomials also depend on four parameters satisfying certain constraint and symmetry that could be used to eliminate one of them, so we will also keep all four parameters $\mu$, $j_1$, $j_2$, and $j_3$. Secondly, the representation depend only on the roots $\Lambda_{j_1}$, $\Lambda_{j_2}$, and $\Lambda_{j_3}$, not on the numbers $j_i$ nor $q^{j_i}$. So, we can fix $j_3=r$ instead of $\Lambda_{j_3}=\Lambda_r$. Finally, it holds that $\Lambda_j=\Lambda_k$ if and only if $q^j=q^k$ or $q^j=q^{2\mu-k+1}$, so the condition $\Lambda_{j_i}\neq\Lambda_r$ can be reformulated as $q^{j_i}\neq q^r$ and $q^{j_i}\neq q^{2\mu-r+1}$.
\smallskip

Our goal is to show the correspondence with orthogonal polynomials. For simplicity, we will work only with the solution \ref[eq.A1jk], \ref[eq.A2jk] and we will not express the form of $I_2$ in the following text. It could be computed easily in a similar way, but we will not need it.
\label[eq.I1vobfin]
$$\eqalignno{I_1v_j=&-q^{\mu-j+j_1+j_2+j_3-1/2}{1+q^{-2\mu+2j}\over (q-q^{-1})^2}A_{j-1}C_jv_{j-1}&\eqmark\cr
&+{\im q^{-(2\mu+2-j_1-j_2-j_3)/2}\over (q-q^{-1})}(A_j+C_j-1+q^{2\mu+2-j_1-j_2-j_3})v_j-{q^{-\mu+j-1/2}\over 1+q^{-2\mu+2j}}v_{j+1},}$$
where
\label[eq.Aj]$$ A_j={(1-q^{j-j_1+1})(1-q^{j-j_2+1})(1-q^{j-j_3+1})(1+q^{-2\mu+j})\over(1+q^{-2\mu+2j})(1+q^{-2\mu+2j+1})},\eqmark$$
\label[eq.Cj]$$C_j={-q^{2\mu+2-j_1-j_2-j_3}(1-q^j)(1+q^{-2\mu-1+j+j_1})(1+q^{-2\mu-1+j+j_2})(1+q^{-2\mu-1+j+j_3})\over(1+q^{-2\mu+2j-1})(1+q^{-2\mu+2j})}.\eqmark$$

\label[secc.00]
\secc Representations for $A_1=A_2=0$

\label[V.AWneklbaze]
\theorem Let $q$ is not a root of unity. Let $R$ be an irreducible representation of $\AW_q(0,0,A_3)$ on $V$, $\dim V=N+1$ satisfying $\ker(I_3-2\epsilon/(q-q^{-1}))=\{0\}$ for $\epsilon=\pm 1$. Then there exists a complex number $\mu$ and a non-zero vector $v_0\in\ker(I_3+\im[\mu]_q)$ such that $O_\mu v_0=0$ and if we define
$$v_{j+1}:=R_{\mu-j}v_j\quad\hbox{for $j=0,1,\dots N$},$$
then the tuple $(v_0,\dots,v_N)$ forms a basis of $V$.
\proof The construction can be performed in the same way as in Theorem \ref[V.AWbaze] (now, we do not have to ensure that $[\mu-j]_q\neq[\mu-j-1]_q$).\qed

The case $A_1=A_2=0$ is, of course, more similar to $U'_q(\so_3)$, where all the parameters are zero. Here, the non-zero parameter $A_3$ causes shift of spectra of the representations.

Solving the system \ref[eq.I1I2soust] we get
\label[eq.I1v00]$$I_1v_j={-q^{-1/2}\over q^{-\mu+j}+q^{\mu-j}}(\tilde D_jv_{j-1}+v_{j+1}),\eqmark$$
\label[eq.I2v00]$$I_2v_j={\im\over q^{-\mu+j}+q^{\mu-j}}(q^{\mu-j}\tilde D_jv_{j-1}-q^{-\mu+j}v_{j+1}).\eqmark$$

Here, we make use of the fact that $\tilde D_j$ itself is a polynomial of degree two in $\Lambda_j$ and factorize it. We get
\label[eq.Dj00]
$$\tilde D_j={q\over(q^{1/2}+q^{-1/2})^2}(\Lambda_j-\Lambda_0)(\Lambda_j+\Lambda_0+\im(q^{1/2}+q^{-1/2})^2A_3).\eqmark$$

Now, consider a classical representation, we will come back to the non-classical case later. We can, of course, use Lemma \ref[L.klasvlč] also in this particular case, so we have $v_{N+1}=0$ and $\tilde D_{N+1}=0$. This holds if $\Lambda_{N+1}=\Lambda_0$ or $\Lambda_{N+1}=-\Lambda_0-\im(q^{1/2}+q^{-1/2})^2A_3$. The first possibility was already excluded in the general case.

So, assume $\Lambda_{N+1}=-\Lambda_0-\im(q^{1/2}+q^{-1/2})^2A_3$. Rearranging the equality we get
\label[eq.nyA3]
$$(q^{-(N+1)/2}+q^{(N+1)/2})[\mu-N/2]_q=-\im(q^{1/2}+q^{-1/2})A_3.\eqmark$$
The bracket on the left-hand side cannot be zero since $q$ is not a root of unity. Thus, it is a quadratic equation in $q^{\mu-N/2}$. In Lemma \ref[L.Rekviv] we have already proven that the representations corresponding to those two solutions have to be equivalent.

Finally, we have to decide when the representation is irreducible.
In Lemma \ref[L.klasvlč] we have shown that in the classical case we have $q^{2\mu}\neq -q^l$, where $l\in\{-1,0,1,\dots,2N+1\}$ and the operator $I_3$ has mutually different eigenvalues. We can show the irreducibility in the same way as in the general case. The condition we mentioned can be rewritten in terms of the parameter $A_3$ by substituting into the equation \ref[eq.nyA3]:
\label[eq.A3ner]
$${(q^{(N+1)/2}+q^{-(N+1)/2})(q^{(k-N)/2}+q^{-(k-N)/2})\over q-q^{-1}}\neq -\epsilon(q^{1/2}+q^{-1/2})A_3,\quad k\in\{-1,0,1,\dots,2N+1\}.\eqmark$$

The final formulas for the representation can be expressed in the following form:
\label[eq.I1vklas00]$$\eqalign{I_1v_j=&{-q^{\mu-j+3/2}\over (1+q^{-2\mu+2j})(q-q^{-1})}(1-q^j)(1-q^{j-N-1})\cr&(1+q^{-2\mu+j-1})(1+q^{-2\mu+j+N})v_{j-1}+{-q^{-\mu+j-1/2}\over 1+q^{-2\mu+2j}}v_{j+1},}\eqmark$$
\label[eq.I2vklas00]$$\eqalign{I_2v_j=&{-\im q^{2\mu-2j+2}\over (1+q^{-2\mu+2j})(q-q^{-1})}(1-q^j)(1-q^{j-N-1})\cr&(1+q^{-2\mu+j-1})(1+q^{-2\mu+j+N})v_{j-1}+{-\im q^{-2\mu}\over 1+q^{-2\mu+2j}}v_{j+1},}\eqmark$$
\label[eq.I3vklas00]$$I_3v_j=-\im{q^{\mu-j}-q^{-\mu+j}\over q-q^{-1}}v_j,\eqmark$$
where $\mu$ is an arbitrary number satisfying \ref[eq.nyA3], $j\in\{0,\dots,N\}$, and $x_{-1}=x_{N+1}=0$. We can see that those representations correspond to the general ones analysed in the previous section for $q^{j_1}=q^r$, $q^{j_2}=-\im q^{\mu+1/2}$, and $q^{j_3}=\im q^{\mu+1/2}$ (the order is, of course, irrelevant).

Now, we move to the case of non-classical representations.

\lemma Consider the same notation as in Theorem \ref[V.AWneklbaze]. Let $R$ be an $(N+1)$-dimensional non-classical irreducible representation of $\AW_q(0,0,A_3)$ such that $\ker(I_3-2\epsilon/(q-q^{-1}))=\{0\}$ for $\epsilon=\pm 1$. Than the operator $I_3$ is diagonalizable and has mutually different eigenvalues, which correspond to eigenvectors $v_0,\dots,v_N$. In addition, we have $[\mu-N]_q=[\mu-N-1]_q$, i.e. $q^\mu=\im\epsilon q^{N+1/2}$ for $\epsilon\in\{-1,1\}$ and $v_{N+1}=av_N$, where
\label[eq.a]
$$a^2=-q[N+1]_q^2-\epsilon q{(q^{(N+1)/2}-q^{-(N+1)/2})^2\over q-q^{-1}}A_3.\eqmark$$
\proof
From the assumptions it follows that $\epsilon/(q^{1/2}-q^{-1/2})$ is an eigenvalue of $I_3$, so there exists $k\in\{0,\dots,N\}$ such that $-\im[\mu-k]_q=\epsilon/(q^{1/2}-q^{-1/2})$, so $q^\mu=\im\epsilon q^{k\pm 1/2}$, and so $[\mu-k]_q=[\mu-(k\pm 1)]_q$. Equivalently, there exists $k\in\{0,\dots,N+1\}$ such that $q^\mu=\im\epsilon q^{k-1/2}$, i.e. $[\mu-k]_q=[\mu-k+1]_q=\im\epsilon/(q^{1/2}-q^{-1/2})$. From the way of construction of the basis $\{v_j\}$ it is clear that we cannot have $[\mu]_q=[\mu+1]_q$, so we have $k>0$. Next, we show that for $k<N+1$ we have a reducible representation. According to Lemma \ref[L.qnum] the equality $k=N+1$ means that the numbers $[\mu]_q,\dots,[\mu-N]_q$ and hence the eigenvalues of $I_3$ are mutually different.

Consider $q^\mu=\im\epsilon q^{k-1/2}$, $k\in\{1,\dots,N+1\}$. According to Lemma \ref[L.qnum] we have $[\mu-N-1]_q=[\mu-l]_q$ if and only if $l=N+1$ or $l=2k-N$. Therefore, the vector $v_{N+1}$ lies in an eigenspace corresponding to the eigenvalue $-\im[\mu-2k+N]_q$. This subspace is one-dimensional for $k>(N+1)/2$, otherwise it is trivial. Thus, we can write $v_{N+1}=av_{2k-N}$, where $v_j=0$ for $j\le 0$ and, for simplicity, we choose $a=0$ in this case, otherwise we have $a\in\C$.

Substituting \ref[eq.I1v00], \ref[eq.I2v00], and \ref[eq.I3v] into condition $(q^{1/2}I_1I_2-q^{-1/2}I_2I_1-I_3-A_3)v_N$ (relation \ref[eq.AW1]) we get
\label[eq.Dcond00nekl]
$${\tilde D_{N+1}\over q^{N-k+1/2}+q^{-N+k-1/2}}v_N+a{1\over q^{N-k+1/2}-q^{-N+k-1/2}}v_{2k-N-1}=0.\eqmark$$
This is satisfied if $a=0$ or $2k-N-1=N+1$. The equality $2k-N-1=N+1$, i.e. $k=N+1$ means $v_{2k-N-1}=v_{N+1}=av_N$ and $q^\mu=\im\epsilon q^{N+1/2}$. Substituting in \ref[eq.Dcond00nekl] we get $a^2=-\tilde D_{N+1}$. Substituting in \ref[eq.Dj00] we get the equality \ref[eq.a].

Now, consider $k<N+1$ and $a=0$, so $k\in\{1,\dots,N\}$ and $v_k=0$, we show that the representation is reducible. According to Lemma \ref[L.qnum] we have $[\mu-j]_q=[\mu-2k+j+1]_q$ for all $j\in{\bbchar Z}$. From the preceding equality it is easy to show that $\Lambda_j=\Lambda_{2k-j}$, so $\tilde D_j=\tilde D_{2k-j}$. Define
$$w_j=\prod_{l=1}^j\tilde D_l^{-1/2}\cdot v_j+\prod_{l=1}^{2k-j-1}\tilde D_l^{-1/2}\cdot v_{2k-j-1},$$
for $j\in{\bbchar Z}$, where an empty product is equal to one and $v_j=0$ for $j\not\in\{0,\dots,N\}$. It is a linear combination of vectors corresponding to the same eigenvalue (or zero vectors), so
$$I_3w_j=-\im[\mu-j]_qw_j.$$
We can also express
$$I_1w_j={-\im\epsilon q^{-1/2}\over q^{-k+j+1/2}+q^{k-j-1/2}}(\tilde D_j^{1/2}w_{j-1}+\tilde D_{j+1}^{1/2}w_{j+1}),$$
$$I_2w_j={-\im\over q^{-k+j+1/2}-q^{k-j-1/2}}(q^{k-j-1/2}\tilde D_j^{1/2}w_{j-1}+q^{-k+j+1/2}\tilde D_{j+1}^{1/2}w_{j+1}).$$
Since we have also $w_j=w_{2k-j-1}$ and, in particular, $w_{k-1}=w_{k}$ and also $w_{-1}=0$ for $k\ge (N+1)/2$ or $w_{-N-1}=0$ for $k\le (N+1)/2$, we can see that the span of the vectors $w_{k-1},\dots, w_0$ for $k\ge N$ or $w_{k},\dots,w_N$ for $k\le N$ forms an invariant subspace of the representation.
\qed

Since $a$ is determined by \ref[eq.a] up to sign, we found additional four representations of the algebra. The proof of the irreducibility is the same as in the case of classical representations. The representations with different $\epsilon$ have different spectra of $I_3$. The representations with different $a$ have different traces of $I_1$. Just in case of $A_3=\mp(q^{(N+1)/2}+q^{-(N+1)/2})^2$ and $\epsilon=\pm 1$ we have $a=0$, so there are only three non-classical representations.

To write the final formulas we have to express
$$\eqalignno{\tilde D_j&=q{(q^{j/2}-q^{-j/2})(q^{N+1-j/2}-q^{-N-1+j/2})\over q-q^{-1}}\left({(q^{j/2}+q^{-j/2})(q^{N+1-j/2}+q^{-N-1+j/2})\over q-q^{-1}}+\epsilon A_3\right)\cr
&=q[j]_q[2N+2-j]_q+q\epsilon{(q^{j/2}-q^{-j/2})(q^{N+1-j/2}-q^{-N-1+j/2})\over q-q^{-1}}A_3.&\eqmark}$$
We are going to use the first row, but the second form illustrates the transition to the representations of $U'_q(\so_3)$ (i.e. for $A_3=0$) listed in \cite[Uqso3].
\label[eq.I1vneklas00]$$\eqalign{I_1v_j=&{\im\epsilon\over (1-q^{-2N+2j-1})}\biggl(-q{(1-q^j)(1-q^{-2N-2+j})\over q-q^{-1}}\cr&\left({q^{N+1-j}(1+q^{j})(1+q^{-2N-2+j})\over q-q^{-1}}+\epsilon A_3\right)v_{j-1}+q^{-N-1+j}v_{j+1}\biggr),}\eqmark$$
\label[eq.I2vneklas00]$$\eqalign{I_2v_j=&{\im\over (1-q^{-2N+2j-1})}\biggl(-q^{N-j+2}{(1-q^j)(1-q^{-2N-2+j})\over q-q^{-1}}\cr&\left({q^{N+1-j}(1+q^{j})(1+q^{-2N-2+j})\over q-q^{-1}}+\epsilon A_3\right)v_{j-1}-q^{-2N+2j-1}v_{j+1}\biggr),}\eqmark$$
\label[eq.I3vneklas00]$$I_3v_j=\epsilon{q^{N-j+1/2}-q^{-N+j-1/2}\over q-q^{-1}}v_j,\eqmark$$
where $\epsilon\in\{-1,1\}$, $j\in\{0,\dots,N\}$, $v_{-1}=0$, and $v_{N+1}=av_N$, where $a$ satisfies \ref[eq.a].

Note that the form of the representation coincide with the classical representations with $q^\mu=\im\epsilon q^{N+1/2}$, $j_1=j_2=r$ and $j_3$ determined by the equation $\Lambda_{j_3}=\Lambda_0+\im(q^{1/2}+q^{-1/2})^2A_3$. The only differences are following. The equation \ref[eq.Djjk] or the equations \ref[eq.Aj], \ref[eq.Cj] contain an removable singularity for $j=N$ that was caused by expanding the formula for $\tilde D_j$ by $g_j^2$. Second difference with the classical representations is the fact that we have $\tilde D_{N+1}\neq 0$ and $v_{N+1}\neq 0$, which causes an extra term in matrix representation of $I_1$ and $I_2$. In the case when $A_1=A_2=0$ the matrices are tridiagonal with zero diagonal except the very last entry.

This completes the classification of all representations of $\AW(0,0,A_3)$ for which $2\epsilon/(q-q^{-1})$ is not an eigenvector of $I_3$.

\theorem Let $q$ is not a root of unity. Let $R$ be an irreducible representations of $\AW_q(0,0,A_3)$ such that $2\epsilon/(q-q^{-1})$ is not an eigenvalue of $I_3$. Assuming inequality \ref[eq.A3ner] the representation is equivalent to one of the five non-equivalent representations---the classical one given by equations \ref[eq.I1vklas00]--\ref[eq.I3vklas00] or one of the four non-classical ones given by equations \ref[eq.I1vneklas00]--\ref[eq.I3vneklas00]. Assuming equality in the relation \ref[eq.A3ner] for some $k$, the representation has to be equivalent to one of the four (or three if $k=2N+1$) non-classical ones given by equations \ref[eq.I1vneklas00]--\ref[eq.I3vneklas00].
\smallskip

\label[sec.ekvivrep]
\sec Dual representations

First of all, we define a representation equivalent to the representation constructed in Section \ref[secc.obrep] multiplying the basis vectors by some scalar. The result will be a bit more symmetric. We define a basis $\X=\{x_k\}_{k=0}^N$ as
$$x_j=(1+q^{-2\mu+2j})\prod_{k=0}^j\im q^{-k+(j_1+j_2+j_3+1)/2}A_{k-1}{1\over (q-q^{-1})}\cdot v_j.$$
The equation \ref[eq.I1vobfin] will change to
\label[eq.I1xobfin]
$$I_1x_j=-{\im q^\nu\over q-q^{-1}}C_jx_{j-1}+{\im q^\nu\over q-q^{-1}}(A_j+C_j-1+q^{-2\nu})x_j-{\im q^{\nu}\over q-q^{-1}}A_jx_{j+1},\eqmark$$
where $\nu=-\mu-1+(j_1+j_2+j_3)/2$. The action of $I_3$ will, of course, not change, so
\label[eq.I3xobfin]
$$I_3x_j=-\im[\mu-j]_qx_j.\eqmark$$

Now, we will try to construct an equivalent representation that would define a Leonard pair. Consider a representation of algebra $\AW_q(A_1,A_2,A_3)$ defined by equations \ref[eq.I1xobfin], \ref[eq.I3xobfin] and numbers $\mu$, $j_1$, $j_2$, and $j_3$. Our goal is to find an equivalent representation, where $I_1$ is diagonal and $I_3$ irreducible tridiagonal. We will make use of the representation of algebra $\AW_q(A_3,A_2,A_1)$ defined by numbers $\nu=-\mu-1+(j_1+j_2+j_3)/2$, $j_1$, $j_2$, $j_3$. Substituting into \ref[eq.A3jk]--\ref[eq.A2jk] we can check that it is indeed a representation of algebra with parameters $A_3,A_2,A_1$. Using isomorphism $\tilde\sigma=\rho^{-1}\sigma\rho$ mapping $I_1\mapsto I_3$ and $I_3\mapsto I_1$ we finally get a representation of $\AW_q(A_1,A_2,A_3)$ we are looking for.

\label[L.dual]
\lemma Let $\mu$, $j_1$, $j_2$, $j_3$, $\nu:=-\mu-1+(j_1+j_2+j_3)/2$ be complex numbers satisfying the following conditions: $q^{2\mu}\neq-q^l$ for every $l\in\{-1,0,1,\dots,2N+1\}$, one of the numbers $q^{j_1}$, $q^{j_2}$, $q^{j_3}$ is equal to $q^{N+1}$ and none of the numbers $\Lambda_{j_1}$, $\Lambda_{j_2}$, nor $\Lambda_{j_3}$ is equal to $\Lambda_k$ for every $k\in\{1,\dots,N\}$. Suppose the same inequalities hold after $\mu$ and $\nu$ are interchanged. Then the classical irreducible representation of $\AW_q(A_1,A_2,A_3)$, where the parameters $A_1$, $A_2$, and $A_3$ are defined by equations \ref[eq.A3jk]--\ref[eq.A2jk], derived in Section \ref[secc.obrep], is equivalent to a representation defined by the same equations after changing
$$I_1\mapsto I_3,\quad I_2\mapsto I_2+{I_1I_3-I_3I_1\over q^{1/2}+q^{-1/2}},\quad I_3\mapsto I_1,\quad \mu\mapsto\nu,\quad \nu\mapsto\mu.$$
Such a representation will be called {\em dual} to the original one.
\proof It is clear that it does not matter if we use the form \ref[eq.I1vobfin] or \ref[eq.I1xobfin], \ref[eq.I3xobfin]. Let us consider the more symmetric form defined above. Then the dual representation has the following form:
\label[eq.I1y]$$I_1y_j=-\im[\nu-j]_qy_j,\eqmark$$
\label[eq.I3y]$$I_3y_j=-{\im q^\mu\over q-q^{-1}}D_jy_{j-1}+{\im q^\mu\over q-q^{-1}}(B_j+D_j-1+q^{-2\mu})y_j-{\im q^{\mu}\over q-q^{-1}}B_jy_{j+1},\eqmark$$
where
\label[eq.Bj]$$ B_j={(1-q^{j-j_1+1})(1-q^{j-j_2+1})(1-q^{j-j_3+1})(1+q^{-2\nu+j})\over(1+q^{-2\nu+2j})(1+q^{-2\nu+2j+1})},\eqmark$$
\label[eq.Dj]$$D_j={-q^{-2\mu}(1-q^j)(1+q^{-2\nu-1+j+j_1})(1+q^{-2\nu-1+j+j_2})(1+q^{-2\nu-1+j+j_3})\over(1+q^{-2\nu+2j-1})(1+q^{-2\nu+2j})}.\eqmark$$

As we mentioned, by interchanging $\nu$ and $\mu$, we get an irreducible representation of algebra $\AW_q(A_3,A_2,A_1)$. Thus, applying the isomorphism $\tilde\sigma$, we indeed obtain an irreducible representation of $\AW_q(A_1,A_2,A_3)$. If it is classical, then, according to Theorem \ref[V.obecklas], it has to be equivalent to one of the representations we have already found. These are, according to Lemma \ref[L.Rekviv], determined by trace. Hence, we only have to show that the dual representation is classical and that $I_3$ has the same trace as the original one. By means of direct computation we can check that the trace of $I_3$ in dual representation is indeed $-\im(q^{(N+1)/2}-q^{-(N+1)/2})/(q^{1/2}-q^{-1/2})\;[\mu-N/2]_q$. Now, we show that $I_3$ has the same eigenvalues in dual representation as in the original one, namely $-\im[\mu]_q,\dots,-\im[\mu-N]_q$.

Firstly, we show that $-\im[\mu]_q$ is an eigenvalue of $I_3$, which means $\det(I_3+\im[\mu]_q)=0$. Using the form of the representation \ref[eq.I3y] only, we can show using induction on the dimension $r$ that in the dual representation we have
$$\det(I_3+\im[\mu]_q)={-\im q^\mu\over q-q^{-1}}{\prod_{k=1}^3 (q^{1-j_k}-1)(q^{2-j_k}-1)\cdots(q^{r-j_k}-1)\over(1+q^{-2\mu+N+1})(1+q^{-2\mu+N+2})\cdots(1+q^{-2\mu+2N+1})}.$$
Using the assumption that $q^{j_k}=q^{N+1}$ for some $k$ we get the zero.

Without loss of generality, we can again assume that the numbers $[\mu-N+1]_q$, $[\mu-N]_q$, \dots, $[\mu]_q$, $[\mu+1]_q$, \dots are all mutually different. We can, therefore, repeat the construction of the eigenbasis for the dual representation as well. Let $\tilde w_0$ be an eigenvector corresponding to the eigenvalue $-\im[\mu]_q$. We will apply $O_{\mu+j}$ repeatedly until we get a linearly dependent vector. The last linearly independent vector will be denoted $\tilde v_0$ and the corresponding eigenvalue $-\im[\tilde\mu]_q$. Then we define $\tilde v_{j+1}:=R_{\tilde\mu-j}\tilde v_j$. As we mentioned, we cannot have $[\tilde\mu-j]_q=[\tilde\mu-j-1]_q$, so we have again constructed an eigenbasis $\{\tilde v_0,\dots,\tilde v_N\}$. The corresponding eigenvalues are $-\im[\tilde\mu]_q,\dots,-\im[\tilde\mu-N]_q$. Computing their sum, we get the trace of $I_3$ in dual representation. Comparing with the previous computation we get $\tilde\mu=\mu$.
\qed

\label[P.superneklas]
\remark To construct a dual representation, we do not need to assume the existence of the original one. It is sufficient to fulfil assumptions for the existence of classical representation defined by numbers $\nu$, $j_1$, $j_2$, and $j_3$, we do not have to exclude the possibility of $q^\mu=\im\epsilon q^{l/2}$ for some $l\in\{0,\dots,2N+1\}$. In that case, an eigenbasis for $I_3$ does not have to exist. Nevertheless, it still holds that the numbers $[\mu-\lfloor l/2\rfloor]_q,[\mu-\lfloor l/2\rfloor+1]_q,\dots$ are mutually different, so choosing $\tilde v_0$ an eigenvector corresponding to the eigenvalue $-\im[\mu]_q$ we can define $\tilde v_{j+1}=R_{\mu-j}\tilde v_j$ until we get $0\neq\tilde v_{\lfloor l/2\rfloor}$, which is for $l$ even an eigenvector corresponding to the eigenvalue $-\im[\mu-l/2]_q=2\epsilon/(q-q^{-1})$ and for $l$ odd it is an eigenvector corresponding to $-\im[\mu-(l-1)/2]_q=\epsilon/(q^{1/2}-q^{-1/2})$. (The equality $v_{\lfloor l/2\rfloor}=0$ would contradict the irreducibility of representation.)

This example shows that for certain parameters of Askey--Wilson algebra there exist representations containing both $2\epsilon/(q-q^{-1})$ or $\epsilon/(q^{1/2}-q^{-1/2}$ as eigenvalues of $I_3$.

\sec Correspondence with $q$-Racah polynomials

Now we are going to show that a representation of Askey--Wilson algebra together with its dual representation defines a Leonard pair corresponding  to $q$-Racah polynomials. More detailed study of this Leonard pair is available in \cite[TerwilligerRacah].

Firstly, let us recall the explicit formula for $q$-Racah polynomials that were discovered by Askey and Wilson in \cite[AskeyWilson]. We use the notation from \cite[Koekoek], where properties of all orthogonal polynomial series of the Askey--Wilson scheme are summarized.

\label[eq.racdef]
$$R_n(\mu(x);\alpha,\beta,\gamma,\delta\mid q)=\Fi43(q^{-n},\alpha\beta q^{n+1},q^{-x},\gamma\delta q^{x+1};\alpha q,\beta\delta q,\gamma q|q;q)\quad n=0,1,2,\dots,N,\eqmark$$
where
$$\mu(x)=q^{-x}+\gamma\delta q^{x+1}\quad\hbox{and}\eqmark$$
\label[eq.rackondim]$$\alpha q=q^{-N}\quad\hbox{or}\quad \beta\delta q=q^{-N}\quad\hbox{or}\quad\gamma q=q^{-N},\quad\hbox{where $N\in\N_0$}.\eqmark$$

\lemma Let $\mu$, $j_1$, $j_2$, $j_3$, and $\nu:=-\mu-1+(j_1+j_2+j_3)/2$ be numbers satisfying the same assumptions as in the preceding lemma. Define linear operators $A,B$ on $V$ as
$$Ax_k=C_kx_{k-1}-(A_k+C_k-1+q^{-2\nu})x_k+A_kx_{k+1},\eqmark$$
$$Bx_k=(q^{-k}-q^{-2\mu+k})x_k.\eqmark$$
Then there exists a basis $\Y$, where operators $A$ and $B$ have the following form:
$$Ay_k=(q^{-k}-q^{-2\nu+k})y_k,\eqmark$$
$$By_k=D_ky_{k-1}-(B_k+D_k-1+q^{-2\mu})y_k+B_ky_{k+1}.\eqmark$$
Therefore, the operators $A$ and $B$ form a Leonard pair.
\proof Consider a representation of $\AW$ in the form \ref[eq.I1xobfin], \ref[eq.I3xobfin]. Then we can write
\label[eq.paircorr]$$A=\im q^{-\nu}(q-q^{-1})I_1,\quad B=\im q^{-\mu}(q-q^{-1})I_3.\eqmark$$
The basis $\Y$ correspond to dual representation constructed in Lemma \ref[L.dual].
\qed

Note that this result agrees with \cite[TerVidunas], Theorem 6.2, which essentially says that irreducible representations of Askey--Wilson algebra define a Leonard pair if both the operators have mutually different eigenvalues.

Now, we can show the correspondence to $q$-Racah polynomials. Denote
\label[eq.koresp1]$$\alpha=q^{-j_1},\quad\beta=-q^{-2\mu-1+j_1},\quad\gamma=q^{-j_3},\quad\delta=-q^{2\mu+1-j_1-j_2}.\eqmark$$
Then the $q$-Racah polynomials $R_n(x)$, $n=0,\dots,N$ with parameters $\alpha$, $\beta$, $\gamma$, $\delta$ are hidden in this Leonard pair in the following way
$$\P_{jk}=r_jR_j(\mu(k)),\quad r_j=\prod_{k=0}^j{A_{k-1}\over C_k},$$
where $\P$ is the transition matrix from basis $\Y$ to basis $\X$ and $\mu(x)=q^{-x}+\gamma\delta q^{x+1}$.

Indeed, the sequences $A_n$, $B_n$, $C_n$, and $D_n$ can be expressed in terms of $\alpha$, $\beta$, $\gamma$, and $\delta$ as
$$\eqalignno{
\label[eq.racAn]	A_n&={(1-\alpha q^{n+1})(1-\alpha\beta q^{n+1})(1-\beta\delta q^{n+1})(1-\gamma q^{n+1})\over(1-\alpha\beta q^{2n+1})(1-\alpha\beta q^{2n+2})},&\eqmark\cr
\label[eq.racBn]	B_n&={(1-\alpha q^{n+1})(1-\beta\delta q^{n+1})(1-\gamma q^{n+1})(1-\gamma\delta q^{n+1})\over(1-\gamma\delta q^{2n+1})(1-\gamma\delta q^{2n+2})},&\eqmark\cr
\label[eq.racCn]	C_n&={q(1-q^n)(1-\beta q^n)(\gamma-\alpha\beta q^n)(\delta-\alpha q^n)\over(1-\alpha\beta q^{2n})(1-\alpha\beta q^{2n+1})},&\eqmark\cr
\label[eq.racDn]	D_n&={q(1-q^x)(1-\delta q^n)(\beta-\gamma q^n)(\alpha-\gamma\delta q^n)\over(1-\gamma\delta q^{2n})(1-\gamma\delta q^{2n+1})}.&\eqmark
}$$
The similarity relations $A^\X\P=\P A^\Y$ and $B^\X\P=\P B^\Y$, where $A^\X,A^\Y,B^\X,B^\Y$ denote the matrices of $A$ and $B$ in the bases $\X$ and $\Y$, can be expressed in terms of the $q$-Racah polynomials as
$$A_jR_{j+1}(\mu(k))-(A_j+C_j-1-\gamma\delta q)R_j(\mu(k))+C_jR_{j-1}(\mu(k))=(q^{-k}+\gamma\delta q^{k+1})R_j(\mu(k)),$$
$$(q^{-j}+\alpha\beta q^{j+1})R_j(\mu(k))=D_kR_j(\mu(k-1))-(B_k+D_k-1-\alpha\beta q)R_j(\mu(k))+B_kR_j(\mu(k+1)),$$
which are precisely the three-term recurrence and difference equation for the $q$-Racah polynomials (cf. \cite[Koekoek], eqs. (14.2.3) and (14.2.6)). From the three term recurrence we could also easily compute the orthogonality relation.



Now we interpret the assumptions of the Lemma \ref[L.dual] in terms of the orthogonal polynomials sequence. The condition that one of the numbers $q^{j_i}$ equals to $q^{N+1}$ ensuring the finite dimension of the representation can be formulated as
\label[eq.finpol]$$\alpha q=q^{-N}\quad\hbox{or}\quad \beta\delta q=q^{-N}\quad\hbox{or}\quad\gamma q=q^{-N},\eqmark$$
which ensures finiteness of the orthogonal polynomials series (\cite[Koekoek], eq. (14.2.1)). The conditions $q^{j_i}\neq q^k$, $q^{j_i}\neq q^{2\mu-k+1}$ and $q^{j_i}\neq q^{2\nu-k+1}$ for $k\in\{1,\dots,N\}$ ensuring irreducibility of the representation and hence irreducibility of the matrices in the Leonard pair can be formulated as follows
\label[eq.irpol1]$$\alpha\neq q^{-k},\quad\beta\delta\neq q^{-k},\quad\gamma\neq q^{-k},\eqmark$$
\label[eq.irpol2]$$\beta\neq q^{-k},\quad\alpha\neq\delta q^{-k},\quad\alpha\beta\neq\gamma q^{-k},\eqmark$$
\label[eq.irpol3]$$\gamma\delta\neq\alpha q^{-k},\quad\gamma\neq\beta q^{-k},\quad\delta\neq q^{-k}.\eqmark$$
Although those conditions are usually not mentioned in the literature, they are necessary to obtain orthogonal polynomial sequence with respect to quasi-definite moment functional. If they are not satisfied, one of the coefficients $A_n$, $B_n$, $C_n$, or $D_n$ may be zero for certain $n$. See also the orthogonality relation (14.2.2) in \cite[Koekoek].

The last assumption that $q^{2\mu}\neq -q^l$ and $q^{2\nu}\neq -q^l$ for all $l\in\{-1,0,1,\dots,2N+1\}$ was made a priori to ensure that the representation is classical and therefore diagonalizable. In terms of the parameters $\alpha, \beta,\gamma,\delta$ it means
\label[eq.classpol]$$\alpha\beta q\neq q^{-l}\quad\hbox{and}\quad\gamma\delta q\neq q^{-l}.\eqmark$$
Looking at the formulas \ref[eq.racAn]--\ref[eq.racDn] it seems that those conditions are necessary for the Leonard pair to be well defined since otherwise there may be zero in one of the denominators. Nevertheless if $\alpha\beta q=q^{-l}$ or $\gamma\delta q=q^{-l}$ for $l\in\{-1,0,2N,2N+1\}$ the singularity is removable. However, for $l\in\{1,\dots,2N-1\}$ the condition is indeed necessary, which is also usually not explicitly stated in the literature. Note also that, for example, if we had $\alpha\beta q=q^{-l}$ for $l\in\{1,\dots,2N-1\}$, the $n$-th polynomial $R_n$ would not be of degree $n$ and the $N$-tuple would not be linearly independent.

\sec Non-classical representations

From the formulas for the Askey--Wilson polynomials, we can now go backwards and derive the missing non-classical representations. Those shell have the same form as the classical representations except for some changes (cf. Section \ref[secc.00]).

First of all $I_3$ has again the spectrum $[\mu]_q$, $[\mu-1]_q$, \dots, $[\mu-N]_q$. For $q^{2\mu}=-q^l$, $l\in\{-1,0,2N,2N+1\}$ the eigenvalues are pairwise distinct. In the end, we can show that the representation with $l=-1$ is equivalent to the representation with $l=2N+1$ and the representation $l=0$ is equivalent with representation $l=2N$ (as in Lemma \ref[L.Rekviv]). Thus, we will work now only with the cases $l\in\{2N,2N+1\}$.

Take the non-classical case $l=2N+1$ and define $I_1$ by formula \ref[eq.I1xobfin] for $j_1$, $j_2$, $j_3$ satisfying the standard assumptions as in Theorem \ref[V.obecklas] and $q^mu=\pm\im q^{N+1/2}$. Since $C_0=0$, we do not have to determine the value of $x_{-1}$. Nevertheless, we have $A_N\neq 0$, so we have to determine $x_r$. By analogy with Section \ref[secc.00], we can guess that $x_{N+1}=ax_N$. So, we can substitute into \ref[eq.I1xobfin]
$$\im q^{-\nu}(q-q^{-1})I_1x_N=C_Nx_{N-1}-((1-a)A_N+C_N-1+q^{-2\nu})x_N.$$
Such form would lead to three-term recurrence
$$((1-a)A_N+C_N-1-\gamma\delta q)R_N(\mu(k))+C_NR_{N-1}(\mu(k))=(q^{-k}+\gamma\delta q^{k+1})R_N(\mu(k))$$

We can suppose that this should not contradict the standard form of three-term recurrence. Notice that for $\alpha\beta q=q^{-2N-1}$ we have $R_N=R_{N+1}$, so we can rewrite the three-term recurrence as
$$-(C_N-1-\gamma\delta q)R_N(\mu(k))+C_NR_{N-1}(\mu(k))=(q^{-k}+\gamma\delta q^{k+1})R_N(\mu(k)).$$

From this, we can conclude that $x_{N+1}=x_N$. Finally, we can check that the formulas for classical representation also define a non-classical representation for $q^{2\mu}=-q^{2N+1}$ if we define $x_{N+1}=x_N$.

In the case $l=2N$ we have $\alpha\beta q=q^{-2N}$ and by similar reasoning we can conclude that $x_{N+1}=x_{N-1}$.

\sec Correspondence with Huang's classification

As we mentioned in the introduction, a complete classification of so-called universal Askey--Wilson algebra appeared in \cite[Huang]. The generating elements are represented by two bidiagonal and one tridiagonal matrices. From those explicit formulas we have the following.

\lemma For any irreducible representation of the Askey--Wilson algebra there exists $\mu\in\C$ such that spectrum of $I_3$ is $\{-\im[\mu-j]\}_{j=0}^N$. On the other hand, for any $\mu\in\C$ there exists such a representation for suitable parameters.
\smallskip

The paper also gives the criterion for diagonalizability of the generating elements.

\nazlemma(\cite[Huang], Lemma 4.6) Let $R$ be an irreducible representation of $\AW$. Then all eigenspaces of $I_3$ are one-dimensional

\noindent{\bf Corollary}\space(\cite[Huang], Lemma 5.1){\bf .} Let $R$ be an irreducible representation of $\AW$ and denote $\mu$ such that $\{-\im[\mu-j]\}_{j=0}^N$ is the spectrum of $I_3$. Then the following are equivalent.
\begitems \style n
* $I_3$ is diagonalizable,
* the numbers $[\mu]$, $[\mu-1]$, \dots, $[\mu-N]$ are pairwise distinct,
* $q^{2\mu}\neq-q^l$, $l\in\{1,\dots,2N-1\}$.
\enditems

From these propostions we can see that the assuption that the representation is classical (i.e. $2\epsilon/(q-q^{-1})$ and $\epsilon/(q^{1/2}-q^{-1/2})$ are not eigenvalues of $I_3$) is sufficient to ensure diagonalizability of $I_3$, but not necessary.

Therefore, our classification contains all diagonalizable representations of the Askey--Wilson algebra except the ``border cases'', when $q^{2\mu}\in\{q^{-1},q^0,q^{2N},q^{2N+1}\}$, which we discussed in the previous section.

%
%

\sec Conclusion

We classified representations of Askey--Wilson algebra satisfying certain conditions allowing us to use the shift operators to construct an eigenbasis of $I_3$. Representations satisfying similar condition for another generating element can be obtained by applying suitable isomorphism of the Askey--Wilson algebra as we indicated in Section \ref[sec.ekvivrep]. Those representations are very important since they define the Leonard pair connected to the $q$-Racah polynomials. In such way, $q$-Racah polynomials (and other types of orthogonal polynomials in the Askey scheme) can be then obtained “at no cost” and this is the most valuable side-effect of solving the classification problem of Askey--Wilson algebra.


\sec Acknowledgements

This work was supported by the Grant Agency of the Czech Technical University in Prague, grant numbers SGS15/215/OHK4/3T/14 and SGS16/239/OHK4/3T/14.

\sec References

\bib[AskeyWilson]  R. Askey, and J. Wilson. A Set of Orthogonal Polynomials That Generalize the Racah Coefficients or $6 - j$ Symbols. {\it SIAM J. Math. Anal.} 1979, 10 (5), 1008-1016. \hlink[http://dx.doi.org/10.1137/0510092]{doi:10.1137/0510092}.
\bib[Bergman1978]  G. M. Bergman. The diamond lemma for ring theory. {\it Adv. Math.} 1978, 29 (2), 178--218. \hlink[http://dx.doi.org/10.1016/0001-8708(78)90010-5]{doi:10.1016/0001-8708(78)90010-5}.
\bib[GromadaPosta] D. Gromada, and S. Po\v sta. Automorphisms of algebras and orthogonal polynomials. {\it Acta Polytechnica.} 2014, 54 (6), 394--397. \hlink[http://dx.doi.org/10.14311/AP.2014.54.0394]{doi:10.14311/AP.2014.54.0394}.
\bib[Uqso399]      M. Havl\'\i\v cek, A. U. Klimyk, and S. Po\v sta. Representations of the cyclically symmetric $q$-deformed algebra $\so_q(3)$. {\it J. Math. Phys.} 1999, 40 (4), 2135--2161. \hlink[http://dx.doi.org/10.1063/1.532856]{doi:10.1063/1.532856}.
\bib[Uqso3]        M. Havl\' \i\v cek, and S. Po\v sta. On the Classification of Irreducible Finite-Dimen\-sional Representations of $U'_q(\so_3)$ Algebra. {\it J. Math. Phys.} 2001, 42 (1), 472--500. \hlink[http://dx.doi.org/10.1063/1.1328078]{doi:10.1063/1.1328078}.
\bib[centeruqso3]  M. Havl\'\i\v cek, and S. Po\v sta. Center of quantum algebra $U_q'(\so_3)$. {\it J. Math. Phys.} 2011, 52 (4), 043521. \hlink[http://dx.doi.org/10.1063/1.3579992]{doi:10.1063/1.3579992}
\bib[Huang]        H. Huang. Finite-Dimensional Irreducible Modules of the Universal Askey--Wilson Algebra. {\it Comm. Math. Phys.} 2015, 340 (3), 959--984. \hlink[http://dx.doi.org/10.1007/s00220-015-2467-9]{doi:10.1007/s00220-015-2467-9}.
\bib[centeraw]     H. Huang. Center of the universal Askey--Wilson algebra at roots of unity. {\it Nucl. Phys. B}. 2016, 909 260--296. \hlink[http://dx.doi.org/10.1016/j.nuclphysb.2016.05.006]{doi:10.1016/j.nuclphysb.2016.05.006}.
\bib[Koekoek]      R. Koekoek, P. A. Lesky, and R. F. Swarttouw. {\it Hypergeometric Orthogonal Polynomials and Their $q$-Analogues.} Berlin Heidelberg: Springer-Verlag, 2010.
\bib[Leonard]      D. A. Leonard. Orthogonal Polynomials, Duality and Association Schemes. {\it SIAM J. Math. Anal.} 1982, 13 (4), 656--663. \hlink[http://dx.doi.org/10.1137/0513044]{doi:10.1137/0513044}.
\bib[TerwilligerLeonard] P. Terwilliger. Two linear transformations each tridiagonal with respect to an eigenbasis of the other. {\it Linear Algebra Appl.} 2001, 330 (1), 149--203. \hlink[http://dx.doi.org/10.1016/S0024-3795(01)00242-7]{doi:10.1016/S0024-3795(01)00242-7}.
\bib[TerwilligerIntro] P. Terwilliger. Introduction to Leonard pairs. {\it J. Comput. Appl. Math.} 2003, 153 (1--2), 463--475. \hlink[http://dx.doi.org/10.1016/S0377-0427(02)00600-3]{doi:10.1016/S0377-0427(02)00600-3}.
\bib[TerwilligerRacah] P. Terwilliger. Leonard pairs and the $q$-Racah polynomials. {\it Linear Algebra Appl.} 2004, 387 235--276. \hlink[http://dx.doi.org/10.1016/j.laa.2004.02.014]{doi:10.1016/j.laa.2004.02.014}.
\bib[TerwilligerAW] P. Terwilliger. The Universal Askey--Wilson Algebra. {\it SIGMA.} 2011, 7 (069), \hlink[http://dx.doi.org/10.3842/SIGMA.2011.069]{doi:10.3842/SIGMA.2011.069}.
\bib[TerVidunas] P. Terwilliger, and R. Vidunas. Leonard pairs and the Askey--Wilson relations, {\it J. Algebra Appl.} 2014, 03, 41.
\bib[Wiegmann] P.B. Wiegmann, and A.V. Zabrodin, Algebraization of difference eigenvalue equations related to $U_q({\rm sl}_2)$, {\it Nuclear Phys. B.} 1995, 451, 699--724,
cond-mat/9501129.
\bib[Zhedanov] A. S. Zhedanov. ``Hidden symmetry'' of Askey--Wilson polynomials. {\it Theoret. and Math. Phys.} 1991, 89 (2), 1146--1157. \hlink[http://dx.doi.org/10.1007/BF01015906]{doi:10.1007/BF01015906}.

\bye